\documentclass[russian,twoside,12pt,fleqn]{article} 
\evensidemargin=\oddsidemargin
\evensidemargin=\oddsidemargin
\usepackage[dvips]{graphicx}
\usepackage[cp1251]{inputenc}           
\usepackage[russian]{babel}             
\usepackage{amsmath,amsthm,amssymb,newlfont}
\usepackage{array,hhline}
\usepackage{cite}                       
\usepackage{mathrsfs}
\renewcommand{\C}{\ensuremath{\mathcal{C}} }          
\newcommand{\R}  {\ensuremath  {\mathbb R^d} }        
\newcommand{\Ei} {\ensuremath {E^{\,i}} }             

\newtheorem{MyTheorem}{Theorem}[section]
\newtheorem{MyCorollary}[MyTheorem]{Corollary}
\newtheorem{MyLemma}[MyTheorem]{Lemma}

\newtheorem{MyRemark}[MyTheorem]{Remark}
\theoremstyle{definition}
\newtheorem{MyDef}[MyTheorem]{Definition}

\newenvironment{emph_}{\bfseries}{}                         

\usepackage{caption}
\newcounter{MyTable}[section]

\DeclareCaptionLabelFormat{MyLabelFormat}{\sc Table \thesection.#2.\ \ }
\DeclareCaptionFormat{MyFormat}{#1#3}
\DeclareCaptionLabelSeparator{MySeparator}{.\\}
\newenvironment*{MyTable}[1]              
{                                         
  \begin{table}[!h]
    \tabcolsep=1 mm
    \stepcounter{MyTable}
    \caption{#1}
  \end{table}
}

\title{Two-step General Linear Methods for \\ Retarded Functional Differential Equations}
\date{}
\author{\\ Anton Tuzov \footnote{Department of Control systems, Siberian State Aerospace University, Krasnoyarsk, Russia, e-mail: tuzov@sibsau.ru.}\\}

\begin{document}
\maketitle
\renewcommand{\abstractname}{Abstract}
\begin{abstract}
  This paper presents a class of Two-Step General Linear Methods for the numerical solution  of Retarded Functional Differential Equations. 
  Explicit methods up to order five are constructed. To avoid order reduction for mildly stiff problems
  the uniform stage order of the methods is chosen to be close to uniform order.

\end{abstract}
\section{Two-step General Linear Methods for \\Ordinary Differential Equations}
For the numerical approximation of the solution $y(t)$ of a system of Ordinary Differential Equations
\begin{equation}
  \begin{split}
    y'(t) &=f(t,y), \quad t \in [t_0,T],\\
    y(t_0)&=y_0,
  \end{split}
\end{equation}
where $f:\ \mathbb R \times \R \longrightarrow \R, \quad y_0 \in \R$,

\noindent
we consider the class of General Linear Methods\cite{Butcher}
\begin{alignat}{2}\label{GLMs}
  Y^{[n]}_i&=\sum\limits_{j=1}^s a_{ij}\, hF^{[n]}_j+\sum\limits_{j=1}^r u_{ij} y^{[n-1]}_j, & \qquad &i=1,\dots ,s,  \notag\\
  y^{[n]}_i&=\sum\limits_{j=1}^s b_{ij}\, hF^{[n]}_j+\sum\limits_{j=1}^r v_{ij} y^{[n-1]}_j, &        &i=1,\dots ,r,  \\
  F^{[n]}_i&=f(t_{n-1}+c_i h,Y^{[n]}_i),                                                     &        &i=1,\dots ,s,  \notag
\end{alignat}
where $y^{[n-1]}_1,\dots,y^{[n-1]}_r$ --- input vectors, available at step number $n$,\\
      $Y^{[n]}_1,\dots,Y^{[n]}_s$     --- stage values,
      $F^{[n]}_1,\dots,F^{[n]}_s$     --- derivative values,
      $a_{ij},\ u_{ij},\ b_{ij},\ v_{ij}.$ --- coefficients of the method.

Let us restrict ourselves to two-step GLMs and choose $r=s+2$,

\mathindent=1.6em
\begin{alignat*}{3}
  y^{[n-1]}_1&\approx y(t_{n-1}), \quad & y^{[n-1]}_2&\approx y(t_{n-2}), \quad & y^{[n-1]}_{2+i}&\approx h y'(t_{n-2}+c_i h),\ i=1,\dots ,s.\quad \text{Then}\\
  y^{[n-1]}_1&=y_{n-1},                 & y^{[n-1]}_2&=y_{n-2},                 & y^{[n-1]}_{2+i}&=h f(t_{n-2}+c_i h,Y^{[n-1]}_i)=hF^{[n-1]}_i,\ i=1,\dots,s.
\end{alignat*}
\mathindent=2.5em

and \eqref{GLMs} takes the form
\begin{alignat}{2}\label{two_step_GLMs_draft_ODEs}
  Y^{[n]}_i&=h \sum\limits_{j=1}^s a_{ij}\, F^{[n]}_j+\ u_{i1}y_{n-1}+u_{i2}y_{n-2}+h\sum\limits_{j=1}^s u_{i,\,2+j} F^{[n-1]}_j, & \qquad &i=1,\dots ,s,  \notag\\
  y_n      &=h \sum\limits_{j=1}^s b_{1j}\, F^{[n]}_j+\ v_{11}y_{n-1}+v_{12}y_{n-2}+h\sum\limits_{j=1}^s v_{1,\,2+j} F^{[n-1]}_j,                          \\
  F^{[n]}_i&=f(t_{n-1}+c_i h,Y^{[n]}_i),                                                     &        &i=1,\dots ,s.                                       \notag
\end{alignat}
In the construction of GLMs it is assumed that $y^{[n-1]}_i=u_i\, y(t_{n-1})+v_i\, hy'(t_{n-1})+O(h^2)$
and 'preconsistensy conditions' holds
\begin{equation}\label{preconsistency_conditions}
  \begin{split}
    Vu=u,\\
    Uu={\mathbf 1}.
  \end{split}
\end{equation}
For \eqref{two_step_GLMs_draft_ODEs} we have $u_1=1,\ v_1=0,\quad u_2=1,\ v_2=-1,\quad u_{2+i}=0,\ v_{2+i}=1,\ i=1,\dots,s$.

\noindent
It follows from \eqref{preconsistency_conditions} that
\begin{align*}
  u_{i2}&=1-u_{i1},\quad i=1,\dots,s,\\
  v_{12}&=1-v_{11}.
\end{align*}

\noindent
Let us denote
\begin{alignat*}{5}
  K^{[n]}_i&:=F^{[n]}_i  &\quad \widetilde{a}_{ij}&:=u_{i,\,2+j}, & \quad u_i&:=u_{i1}, &\quad \Rightarrow u_{i2}&=1-u_i, &\quad j&=1,\dots,s,\ i=1,\dots,s,\\
        b_j&:=b_{1j},    &         \widetilde{b_j}&:=v_{1,\,2+j}, &         v&:=v_{11}, &      \Rightarrow v_{12}&=1-v,   &      j&=1,\dots,s,\\
\end{alignat*}
then the method \eqref{two_step_GLMs_draft_ODEs} satisfying 'preconsistensy conditions'~\eqref{preconsistency_conditions} takes the form
\begin{alignat}{4}\label{two_step_GLMs_ODEs}
  y_n      &=  \ (1-v) y_{n-2}+v y_{n-1}   & &+h\sum\limits_{j=1}^s \widetilde{b_j}   K^{[n-1]}_j & &+h \sum\limits_{j=1}^s b_j   K^{[n]}_j, & \quad                \notag\\
  K^{[n]}_i&=f(t_{n-1}+c_i h,Y^{[n]}_i),   & &                                                    & &                                        &       &i=1,\dots ,s, \\
  Y^{[n]}_i&=  (1-u_i) y_{n-2}+u_i y_{n-1} & &+h\sum\limits_{j=1}^s \widetilde{a}_{ij}K^{[n-1]}_j & &+h\sum\limits_{j=1}^s a_{ij} K^{[n]}_j, &       &i=1,\dots ,s. \notag
\end{alignat}

\newpage
\section{Two-step General Linear Methods for \\Retarded Functional Differential Equations}
We begin with notations introduced in~\cite{Maset}.

\noindent
Let $r\in[0,+\infty)$, and \C be the space of continuous functions $[-r,0]\longrightarrow \R$,
equipped with the maximum(uniform) norm $\|\phi\|=\max\limits_{\theta \in [-r,0]}|\phi(\theta)|, \quad \phi\in\C$,
\; where $| \cdot |$ is an arbitrary norm on \R.

\noindent
Let $u$ be continuous function $[a-r,b)\longrightarrow \R$, where $a<b$. Then $\forall \, t\in[a,b)$
shift function is well defined by $u_t(\theta)=u(t+\theta), \quad \theta\in[-r,0]$, and $u_t\in\C$.

Let us consider a system of Retarded Functional Differential Equations
\begin{equation}\label{RFDEs_draft}
  \begin{split}
     y'(t)           &=f(t,y_t),     \quad       t      \in [t_0,T],\\
     y_{t_0}(\theta) &=\phi(\theta), \qquad      \theta \in [-r,0].
  \end{split}
\end{equation}
where $(t_0,\phi)\in \Omega,\ f:\ \Omega \longrightarrow \R,\quad \Omega \subset \mathbb R \times \C,\ \Omega\ \text{is open set}$.

\noindent
It is assumed that there exists a unique solution of~\eqref{RFDEs_draft}.

We introduce the class of two-step GLMs for RFDEs on the base of approach proposed in~\cite{Maset}.
We can reformulate the method~\eqref{two_step_GLMs_ODEs} for RFDEs~\eqref{RFDEs_draft} as follows

\mathindent=0em
\begin{alignat}{4}\label{two_step_GLMs}
  \eta^{[n]}(\alpha h) &=(1-v(\alpha))\eta^{[n-1]}(0)                 & &+v(\alpha)\eta^{[n-1]}(h)   & &+h\sum\limits_{j=1}^s \widetilde{b_j}(\alpha)K^{[n-1]}_j    & &+h\sum\limits_{j=1}^s b_j(\alpha)K^{[n]}_j,\notag\\
                       &                                              & &                            & &\ \alpha\in[0,1],\notag\\
  K^{[n]}_i            &=f(t_{n-1}+c_i h,{Y_{\quad c_ih}^{[n]\; i}}),\hspace{-2em} & &                            & &                                                            & &\ i=1,\dots ,s,\notag\\
  Y^{[n]\; i}(\alpha h)&=(1-u_i(\alpha))\eta^{[n-1]}(0)               & &+u_i(\alpha)\eta^{[n-1]}(h) & &+h\sum\limits_{j=1}^s \widetilde{a}_{ij}(\alpha)K^{[n-1]}_j & &+h\sum\limits_{j=1}^s a_{ij}(\alpha)K^{[n]}_j,\hspace{-2em}\\
                       &                                              & &                            & &\ \alpha\in [0,c_i],                                        & &\ i=1,\dots,s,\notag\\
  \eta^{[n]}(\theta)   &=\eta^{[n-1]}_h(\theta),                      & &                            & &\ \theta \in [-r,0],\notag\\
  Y^{[n]\; i}(\theta)  &=\eta^{[n-1]}_h(\theta),                      & &                            & &\ \theta \in [-r,0],                                        & &\ i=1,\dots ,s,\notag
\end{alignat}
\mathindent=2.5em
where
\begin{alignat*}{2}
  &\eta^{[n-1]}:\ [-r,h]\longrightarrow \R, \quad K^{[n-1]}_i \in \R                                          &\quad        &\text{are available as approximations}\\
  &                                                                                                           &             &\text{computed in the step } n-1,\\
  &Y^{[n]\; i}                                                                                                &\text{--- }  &\text{stage functions},\\
  &K^{[n]}_i                                                                                                  &\text{--- }  &\text{stage values},\\
  &u_i(\cdot),\ \widetilde{a}_{ij}(\cdot), a_{ij}(\cdot),\quad v(\cdot),\ \widetilde{b_j}(\cdot),\ b_j(\cdot) &\text{--- }  &\text{coefficients of the method}.\\
\end{alignat*}
\begin{center}
  \includegraphics[scale=0.75]{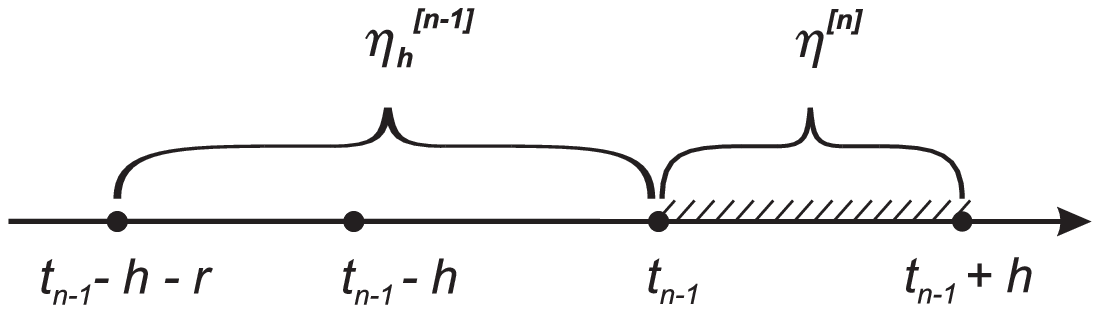}
\end{center}
Thus $\eta^{[n]}(\alpha h) \approx y(t_{n-1}+\alpha h), \ \alpha \in [0,1]$, $\quad \eta^{[n]}(\theta) \approx y(t_{n-1}+\theta),\ \theta \in [-r,0]$, hence\\
$\eta^{[n]}_h\approx y_{t_n}$ on $[-r-h,0]$, where $t_n=t_{n-1}+h$.

\begin{MyRemark}\label{Why_two_step_methods}
  We chose two-step methods among multi-step methods ($k \ge 2$) for the following reasons.
  \begin{itemize}
  \item         
    For multistep methods, the local error $E(h,t_{n-1},y_{t_{n-1}})$ has the required order
    only if exact solution $y(t)$ is sufficiently smooth on $[t_{n-k},t_{n}]$.
    This is rather severe assumption for many problems~\eqref{RFDEs_draft}.
    For the case of two-step methods ($k=2$) this condition imposes the weakest restriction on stepsize.
  \item         
    Furthermore, in the case of two-step methods, starting procedure and stepsize strategy seem to be simplest ones.
  \end{itemize}
\end{MyRemark}

\noindent       
Let us denote
\begin{equation*}
  \eta:=\eta^{[n]},\ K_i:=K^{[n]}_i,\quad \overline{\eta}:=\eta^{[n-1]},\ \overline{K}_i:=K^{[n-1]}_i,\quad \sigma:=t_{n-1},
\end{equation*}
then the method~\eqref{two_step_GLMs} can be reformulated in Stefano Maset's notations as follows

\mathindent=1em
\begin{alignat}{5}\label{TS_GLMs}
  \eta(\alpha h) &=(1-v(\alpha))\overline{\eta}(0) & &+v(\alpha)\overline{\eta}(h)     & &+h\sum\limits_{j=1}^s \widetilde{b_j}(\alpha)\overline{K}_j    & &+h\sum\limits_{j=1}^s b_j(\alpha)K_j,    &\quad &\alpha\in[0,1],\notag\\
  K_i            &=f(\sigma+c_i h,{Y^i_{c_ih}}),   & &                                 & &                                                               & & \notag\\
  Y^i(\alpha h)  &=(1-u_i(\alpha))\overline{\eta}(0) & &+u_i(\alpha)\overline{\eta}(h) & &+h\sum\limits_{j=1}^s \widetilde{a}_{ij}(\alpha)\overline{K}_j & &+h\sum\limits_{j=1}^s a_{ij}(\alpha)K_j, &      &\alpha\in [0,c_i],\\
  \eta(\theta)   &=\overline{\eta}_h(\theta),      & &                                 & &                                                               & &                                         &      &\theta \in [-r,0],\notag\\
  Y^{i}(\theta)  &=\overline{\eta}_h(\theta).      & &                                 & &                                                               & &                                         &      &\theta \in [-r,0],\notag
\end{alignat}
\mathindent=2.5em
\begin{center}
  \includegraphics[scale=0.75]{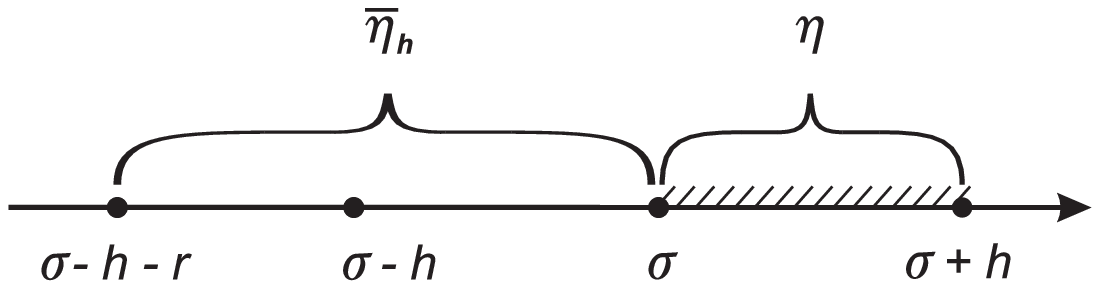}
\end{center}

\newpage
\section{Two-step GLMs for RFDEs in Stefano Maset's notations}
%
%

\noindent
Let us consider a system of Retarded Functional Differential Equations
\begin{equation}\label{RFDEs} 
  \begin{split}
     x'(t)           &=f(t,x_t),     \quad       t      \in [t_0,T],\\
     x_{t_0}(\theta) &=\phi(\theta), \qquad      \theta \in [-r,0].
  \end{split}
\end{equation}
where $(t_0,\phi)\in \Omega,\ f:\ \Omega \longrightarrow \R,\quad \Omega \subset \mathbb R \times \C,\ \Omega\ \text{is open set}$.


It is assumed that conditions of existence and uniqueness theorem for the~\eqref{RFDEs} hold.

\noindent


When $s$-stage Two-Step General Linear Method for RFDEs (TSGLM) with coefficients
$(a_{ij}(\cdot),b_j(\cdot),c_i,\widetilde{a}_{ij}(\cdot),\widetilde{b}_j(\cdot),u_i(\cdot),v(\cdot))_{i,j=1,\dots,s}\ $
is applied with stepsize $h$ to~\eqref{RFDEs} for the computation of the solution $x(t)$,
it yields, as an approximation on $[-r,h]$ of the shift function $y:=x(\sigma+\cdot)$, the function
\mathindent=0em
\begin{alignat}{2}
  \eta(\alpha h) &=(1-v(\alpha))\overline{\eta}(0) +v(\alpha)\overline{\eta}(h)   +h\sum\limits_{j=1}^s \widetilde{b_j}(\alpha)\overline{K}_j    +h\sum\limits_{j=1}^s b_j(\alpha)K_j,&\ \    &\alpha\in[0,1],\\
  \eta(\theta)   &=\overline{\eta}_h(\theta),                                                                                                                                         &       &\theta \in [-r,0],\notag\\
  \notag\\
  \text{where }\notag\\
 &\text{ function } \overline{\eta}_h\approx y \text{ on } [-r-h,0] \text{ and stage values } \overline{K}_i, \text{ are available as approximations}\hspace{-10cm}\notag\\
 &\text{ computed in the previous step,}\notag\\
 K_i             &=f(\sigma+c_i h,{Y^i_{c_ih}}),                                                                                                                                      &       &i=1,\dots ,s,\\
 \notag\\
 &\text{ and } Y^i:\ [-r,c_i h]\longrightarrow \R \ \text{ is a stage function given by}\notag\\
    Y^i(\alpha h)&=(1-u_i(\alpha))\overline{\eta}(0) +u_i(\alpha)\overline{\eta}(h) +h\sum\limits_{j=1}^s \widetilde{a}_{ij}(\alpha)\overline{K}_j +h\sum\limits_{j=1}^s a_{ij}(\alpha)K_j,&      &\alpha\in [0,c_i],\\
    Y^{i}(\theta)&=\overline{\eta}_h(\theta).                                                                                                                                              &      &\theta \in [-r,0],\notag
\end{alignat}
\mathindent=2.5em

It is assumed that coefficients $(a_{ij}(\cdot),b_j(\cdot),c_i,\widetilde{a}_{ij}(\cdot),\widetilde{b}_j(\cdot),u_i(\cdot),v(\cdot))_{i,j=1,\dots,s}\ $
of TSGLMs satisfy the following conditions:
\mathindent=1.5em
\begin{alignat}{4}
  \text{\textbullet}\ & a_{ij}(\cdot),\ \widetilde{a}_{ij}(\cdot),\ u_i(\cdot), &\quad &\text{are polynomial functions } [0,c_i]\longrightarrow {\mathbb R}, & \quad i,\ j&=1,\dots,s.\\
                      & b_j(\cdot),\    \widetilde{b}_j(\cdot),\      v(\cdot), &      &\text{are polynomial functions }   [0,1]\longrightarrow {\mathbb R}, &           j&=1,\dots,s.\notag\\
  \text{\textbullet}\ & c_i \in {\mathbb R},\ c_i \ge 0,                        &      &                                                                     &           i&=1,\dots,s.\\
  \text{\textbullet}\ & a_{ij}(0)=\widetilde{a}_{ij}(0)=0,                      &      & u_i(0)=1,                                                           &       i,\ j&=1,\dots,s.  \label{Yi_continuous}\\
  \text{\textbullet}\ & b_j(0)\ =\widetilde{b}_j(0)\ =0,                        &      & v(0)=1,                                                             &           j&=1,\dots,s.  \label{eta_continuous}
\end{alignat}
\vspace{-12pt}
\mathindent=2.5em

\noindent
The last two conditions correspondingly 
gurantee continuity of the stage functions $Y^i_{c_i h} \in \C$ and
the approximate solution $\eta_h \in \C$ provided that
approximate solution computed in the previous step is continuous function $\overline{\eta}_h \in \C.$

\begin{MyRemark}
  If the conditions
  \begin{alignat}{4}\label{one_step_conditions}
    u_i(\cdot)&=1, & \quad \widetilde{a}_{ij}(\cdot)&=0,\quad & i,\ j&=1,\dots,s,\\
      v(\cdot)&=1, &       \widetilde{b}_j(\cdot)   &=0,\quad &     j&=1,\dots,s,\notag
  \end{alignat}
  hold, the two-step method~\eqref{TS_GLMs} becomes the one-step RK method for RFDEs 
  introduced in~\cite{Maset}, where initial function $\ \phi:=\overline{\eta}_h$.
\end{MyRemark}

\begin{MyDef}
  TSGLM with coefficients $(a_{ij}(\cdot),b_j(\cdot),c_i,\widetilde{a}_{ij}(\cdot),\widetilde{b}_j(\cdot),u_i(\cdot),v(\cdot))_{i,j=1,\dots,s}\ $
  is called explicit if $a_{ij}(\cdot)=0 \text{ for all } j:\ j\ge i,\quad i,j=1,\dots,s.$
\end{MyDef}

\begin{MyDef}
  The function $E=\eta-y:\ [0,h]\longrightarrow \R$ computed under the assumption that
  $\overline{\eta}_h=y$ on $[-r-h,0]$ is called the local error of TSGLM~\eqref{TS_GLMs}.
\end{MyDef}

\begin{MyDef}
  The function $\Ei=Y^i-y:\ [0,c_ih]\longrightarrow \R$ computed under the assumption that
  $\overline{\eta}_h=y$ on $[-r-h,0]$ is called the local stage error of TSGLM~\eqref{TS_GLMs}.
\end{MyDef}

By analogy with~\cite{Dekker} we define \begin{emph_}stage order \end{emph_} for RFDEs.
\begin{MyDef}
  Let functions $\widetilde{E}^{\,i}=Y^i-y:\ [0,c_ih]\longrightarrow \R,\quad i=1,\dots,s$,\\
    $\widetilde{E}^{s+1}=\eta-y:\ [0,h]\longrightarrow \R$ are computed under the assumption that\\
    $\overline{\eta}_h=y$ on $[-r-h,0],\ \overline{K}_j=y'(-h+c_jh),\quad {K}_j=y'(c_jh),\quad j=1,\dots,s,\ $ that is

\mathindent=0em
\begin{alignat}{4}\label{Ei_stage_order}
  \widetilde{E}^{\,i}(\alpha h) &=(1-u_i(\alpha))y(-h) & &+u_i(\alpha)y(0)     & &+h\sum\limits_{j=1}^s \widetilde{a}_{ij}(\alpha)y'(-h+c_j h) & &+h\sum\limits_{j=1}^s a_{ij}(\alpha)y'(c_j h)-\notag\\
                                &-y(\alpha h),         & &                     & &\ \alpha\in [0,c_i],                                         & &\ i=1,\dots,s, \\
  \widetilde{E}^{s+1}(\alpha h) &=(1-v(\alpha))y(-h)   & &+v(\alpha)y(0)       & &+h\sum\limits_{j=1}^s \widetilde{b}_j(\alpha)y'(-h+c_j h)    & &+h\sum\limits_{j=1}^s b_j(\alpha)y'(c_j h)-\notag\\
                                &-y(\alpha h),         & &                     & &\ \alpha\in[0,1] \notag.
\end{alignat}
\mathindent=2.5em
  Denote $c_{s+1}:=1$.
  If there are positive integers $\widetilde{p}_i$ and reals $D_{\,i}>0,\ H>0$ such that
  \begin{equation}
    \max\limits_{\alpha\in[0,c_i]}|\widetilde{E}^{\,i}(\alpha h)|\le D_i\, h^{\widetilde{p}_i+1},\qquad h\in[0,H], \quad i=1,\dots,s+1,
  \end{equation}
\noindent
then positive integer $\widetilde{p}=\min\{\widetilde{p}_1,\dots,\widetilde{p}_{s+1}\}$ is called uniform stage order
of TSGLM~\eqref{TS_GLMs}.
\end{MyDef}

\newpage
\section{Order conditions}\label{sec_Order cond}
Assume that $f$ is of class $C^l$ with respect to the second argument for a sufficiently large $l$
and
solution $x(t)$ of~\eqref{RFDEs} is of piecewise class $C^m$ for a sufficiently large $m$.

We introduce the polynomial functions $\Gamma_k:\ [0,1]\longrightarrow \mathbb{R}$ and $\Gamma_{ik}:\ [0,c_i]\longrightarrow \mathbb{R}$ given by

\mathindent=0em
\begin{alignat}{4}
  \Gamma_{k}(\alpha) &=\frac{1}{(k-1)!}\biggl[\frac{(1-v(\alpha))(-1)^k}{k}   &   &+\sum\limits_{j=1}^{s}{\widetilde{b}_j}(\alpha)(-(1-c_j))^{k-1}    & &+\sum\limits_{j=1}^{s} b_j(\alpha)c_j^{k-1}    & &-\frac{\alpha^k}{k} \biggr],\notag\\
                     &                                                        &   &\ \alpha\in[0,1], \notag\\
  \\
  \Gamma_{ik}(\alpha)&=\frac{1}{(k-1)!}\biggl[\frac{(1-u_i(\alpha))(-1)^k}{k} &   &+\sum\limits_{j=1}^{s}{\widetilde{a}_{ij}}(\alpha)(-(1-c_j))^{k-1} & &+\sum\limits_{j=1}^{s} a_{ij}(\alpha)c_j^{k-1} & &-\frac{\alpha^k}{k} \biggr],\notag\\
                     &                                                        &   &\ \alpha\in[0,c_i],                                                & &\ i=1,\dots,s.\notag
\end{alignat}
\mathindent=2.5em
\begin{MyRemark}
  If the conditions~\eqref{one_step_conditions} hold the $\Gamma_{ik},\ \Gamma_{k}$
  are the same as for the one-step RK method~\cite{Maset}.
\end{MyRemark}

Let $c_1^*,\dots,c_{s^*}^*$ such that $c_1^*<c_2^*<\dots<c_{s^*}^*$ and
$\{c_1^*,\dots,c_{s^*}^*\}=\{c_1,\dots,c_s \}$, i.e. $c_i^*$ are distinct $c_i$ in increasing order.

\begin{MyLemma}\label{local_error_cond}
  Let $p$ be a positive integer. If $x$ is of piecewise class $C^{p+1}$ and
  the local stage errors computed in the previous step are $\overline{E}^{\,i}=O(h^p),\ i=1,\dots,s$, then
  the local error $E$ and the local stage errors $\Ei$ satisfy

  \mathindent=0.5em
  \begin{alignat}{5}
    E(\alpha h)      &=h\sum\limits_{j=1}^{s} b_j(\alpha)D_j     & &+\sum\limits_{k=1}^{p}y^{(k)}(0)h^k\Gamma_k(\alpha)    & &+O(h^{p+1}), & \quad &\alpha\in[0,1],   \\
    \Ei(\alpha h)    &=h\sum\limits_{j=1}^{s} a_{ij}(\alpha)D_j  & &+\sum\limits_{k=1}^{p}y^{(k)}(0)h^k\Gamma_{ik}(\alpha) & &+O(h^{p+1}), &       &\alpha\in[0,c_i], & \quad i&=1,\dots,s,\\
    \text{where } \notag\\
                  D_i&=f(\sigma+c_i h,y_{c_i h}+E^{\,i}_{c_i h})-f(\sigma+c_i h,y_{c_i h}),\hspace{-10cm}  & &              & &             &       &                  &       i&=1,\dots,s.
  \end{alignat}
\end{MyLemma}

*The Lemma~\ref{local_error_cond} has been proved, but its proof is omitted here for brevity.*
\vspace{12pt}

In the following we assume that the TSGLM satisfies the conditions $\Gamma_1=0$ and\\
$\Gamma_{i1}=0,\ i=1,\dots,s,$ that is

\begin{alignat}{6}\label{order_1_cond}
  &v(\alpha)-1   & &+\sum\limits_{j=1}^{s}{\widetilde{b}_j}   (\alpha) & &+\sum\limits_{j=1}^{s} b_j(\alpha)      & &\ =\alpha, &\hspace{6.6em} &\alpha\in[0,1],\notag\\
  \\
  &u_i(\alpha)-1 & &+\sum\limits_{j=1}^{s}{\widetilde{a}_{ij}}(\alpha) & &+\sum\limits_{j=1}^{s} a_{ij}(\alpha)   & &\ =\alpha, &             &\alpha\in[0,c_i],\quad i=1,\dots,s.\notag
\end{alignat}

The above condition is an equivalent form of uniform stage order \begin{emph_}one\end{emph_} condition.

\begin{MyTheorem}\label{order_2_cond}
  A TSGLM satisfying~\eqref{order_1_cond} has 
  uniform order two iff $\ \Gamma_2=0$.
\end{MyTheorem}

*The theorem~\ref{order_2_cond} has been proved, but its proof is omitted here for brevity.*
\vspace{12pt}

\begin{MyTheorem}\label{order_3_cond}
  Let TSGLM satisfy~\eqref{order_1_cond} and has uniform order two.

\noindent
  If $\ \Gamma_3=0\ $ and $\sum\limits_{\substack{i=1\\ c_i=c_m^*}}^s {b_i(\alpha) \Gamma_{i\, 2}(\beta)}=0,\quad \alpha\in[0,1],\quad \beta\in[0,c_m^*],\quad m=1,\dots,s^*$

\noindent
  then the method has uniform order three.
\end{MyTheorem}

*The theorem~\ref{order_3_cond} has been proved, but its proof is omitted here for brevity.*
\vspace{12pt}

\begin{MyTheorem}\label{order_4_cond}
  Let TSGLM satisfy~\eqref{order_1_cond} and has uniform order three.
  \mathindent=0em
\begin{alignat}{5}
  \text{If }\ &\Gamma_4=0, \notag\\
              &\sum\limits_{\substack{i=1 \\ c_i=c_m^*}}^s {b_i(\alpha) \Gamma_{i\, 3}(\beta)}=0,                                                           &\quad \alpha&\in[0,1],&\quad &\beta\in[0,c_m^*], &\quad        &              &\quad   m&=1,\dots,s^*,\\
              &\sum\limits_{\substack{i=1 \\ c_i=c_m^*}}^s {\sum\limits_{\substack{j=1\\ c_j=c_l^*}}^s b_i(\alpha) a_{ij}(\beta) \Gamma_{j\, 2}(\gamma)}=0, &      \alpha&\in[0,1],&      &\beta\in[0,c_m^*], &\quad \gamma &\in[0,c_l^*], &\quad l,m&=1,\dots,s^*.\notag
\end{alignat}

\noindent
  then the method has uniform order four.
\end{MyTheorem}
*The theorem~\ref{order_4_cond} has been proved, but its proof is omitted here for brevity.*
\begin{MyTheorem}\label{stage_order}
  TSGLM has uniform stage order~$\widetilde{p}$ iff\\
  $\Gamma_{ik}=0,\ \Gamma_{k}=0,\quad i=1,\dots,s,\quad k=1,\dots,\widetilde{p}.$
\end{MyTheorem}

\noindent
Proof. Follows by Taylor series expansion of functions $\widetilde{E}^{\,i}$ given by~\eqref{Ei_stage_order}.

\vspace{12pt}

The following results can be obtained as corollary of theorems~\eqref{stage_order} and~\eqref{order_3_cond},~\eqref{order_4_cond}.

\begin{MyCorollary}\label{stage_order_3_cond}
  Let TSGLM has uniform stage order two.

\noindent
  If $\ \Gamma_3=0\ $  then the method has uniform order three.
\end{MyCorollary}

\begin{MyCorollary}\label{stage_order_4_cond}
  Let TSGLM has uniform stage order three.

\noindent
  If $\ \Gamma_4=0$ then the method has uniform order four.
\end{MyCorollary}
The results of Corollary~\eqref{stage_order_3_cond} and~\eqref{stage_order_4_cond} can be easily generalized as follows.

\begin{MyTheorem}\label{order_p_cond}
  Let TSGLM has uniform stage order~$\widetilde{p}$.\\
 It has uniform order~$p=\widetilde{p}+1$ iff $\ \Gamma_{\widetilde{p}+1}=0$.
\end{MyTheorem}

\noindent

*The theorem~\ref{order_p_cond} has been proved, but its proof is omitted here for brevity.*
\vspace{12pt}

\newpage
\section{Construction of explicit two-stage GLMs of\\ uniform stage order four and five}
Consider two-stage explicit TSGLM satisfying~\eqref{order_1_cond}. 
It's Butcher tableau is

\begin{MyTable}{Butcher tableau for 2-stage explicit TSGLMs}\label{tab_2_stage_ETSGLMs}
  \arrayrulewidth=0.5pt
  \begin{tabular}{c|c|c c|c c}
    $c_1$   & $u_1(\alpha)$ & $\widetilde{a}_{11}(\alpha)$  & $\widetilde{a}_{12}(\alpha)$  & 0                 & 0             \\
    $c_2$   & $u_2(\alpha)$ & $\widetilde{a}_{21}(\alpha)$  & $\widetilde{a}_{22}(\alpha)$  & $a_{21}(\alpha)$  & 0             \\

    \hline
          & $v(\alpha)$   & $\widetilde{b}_{1}(\alpha)$   & $\widetilde{b}_{2}(\alpha)$   & $b_1(\alpha)$     & $b_{2}(\alpha)$
  \end{tabular}
  \vspace{12pt}
\end{MyTable}
A natural choice will be to space out  the abscissae $c_i,\ i=1,\dots,s$ uniformly
\smallskip
in the interval~$[0,1]$ so that~\cite{Butcher_Jackiewicz}
\smallskip
$c_1=0,\ c_2=\dfrac{1}{s-1},\dots,\ c_{s-1}=\dfrac{s-2}{s-1},\ c_s=1$. In the case of $s=2$ we have $c_1=0,\ c_2=1$.

Since $c_1=0,$ conditions $\Gamma_{1k}(\alpha)=0,\quad\alpha\in[0,c_1],\ k=1,2,\dots$
reduce to
$\Gamma_{1k}(0)=0,$ $\ k=1,2,\dots$ that follows from~\eqref{Yi_continuous}. It also follows that $u_1(\cdot)=1,\ \widetilde{a}_{11}(\cdot)=0,\ \widetilde{a}_{12}(\cdot)=0$.

For brevity we omit the argument $\alpha$ of the method coefficient functions.
By theorem~\eqref{order_p_cond}, the method has uniform order four and uniform stage order three if\\
 $\Gamma_k=0,\quad k=1,2,3,4\ \text{ and } \Gamma_{2\,k}=0,\quad k=1,2,3,\quad $ that is
\smallskip
\begin{alignat}{5}\label{2_stage_4_order_ETSGLM_cond}
  -(1-v)           &+\widetilde{b}_1     & &+\widetilde{b}_2     & &+b_1    & &+b_2      & &=\alpha,            \notag\\
  \frac{1-v}{2}    &-\widetilde{b}_1     & &                     & &        & &+b_2      & &=\frac{\alpha^2}{2},\notag\\
  -\frac{1-v}{3}   &+\widetilde{b}_1     & &                     & &        & &+b_2      & &=\frac{\alpha^3}{3},\notag\\
   \frac{1-v}{4}   &-\widetilde{b}_1     & &                     & &        & &+b_2      & &=\frac{\alpha^4}{4},\\
  -(1-u_2)         &+\widetilde{a}_{21}  & &+\widetilde{a}_{22}  & &+a_{21} & &          & &=\alpha,            \notag\\
  \frac{1-u_2}{2}  &-\widetilde{a}_{21}  & &                     & &        & &          & &=\frac{\alpha^2}{2},\notag\\
  -\frac{1-u_2}{3} &+\widetilde{a}_{21}  & &                     & &        & &          & &=\frac{\alpha^3}{3},\notag
\end{alignat}
where $\alpha\in[0,1]$.

\noindent
The coefficients are defined by
\begin{alignat}{5}\label{2_stage_4_order_ETSGLM_coeff}
u_{{2}}               &=- \left( 2\,\alpha-1 \right)  \left( \alpha+1 \right) ^{2},\notag\\
v                     &= \left( \alpha-1 \right) ^{2} \left( \alpha+1 \right) ^{2},\notag\\
\widetilde{a}_{2,1}   &={\alpha}^{2} \left( \alpha+1 \right),\notag\\
\widetilde{b}_{1}     &=-\frac{1}{12}\,{\alpha}^{2} \left( \alpha+1 \right)  \left( 5\,\alpha-7 \right),\\
a_{2,1}               &= \alpha\, \left( \alpha+1 \right)^2-\widetilde{a}_{2,2},\notag\\
b_1                   &=-\frac{1}{3}\,\alpha\, \left( 2\,\alpha-3 \right)  \left( \alpha+1 \right) ^{2}-\widetilde{b}_{2},\notag\\
b_2                   &=\frac{1}{12}\,{\alpha}^{2} \left( \alpha+1 \right) ^{2},\notag
\end{alignat}
where $\widetilde{a}_{2,2},\ \widetilde{b}_{2}$ remain free.
The relation $\Gamma_5(1)=\dfrac{4}{15}\ne0$ implies that it is impossible to attain discrete order five.

The uniform order and the uniform stage order can be increased by finding a suitable value for $c_2$.
Assume that $c_1=0,\ c_2  \ne 0$ (in general case~$c_2 \ne 1$).
By theorem~\eqref{order_p_cond}, the method has uniform order five and uniform stage order four if\\
$\Gamma_k=0,\quad k=1,2,3,4,5\ \text{ and }$ $\Gamma_{2\,k}=0,\quad k=1,2,3,4,\quad $ that is
\smallskip
\begin{alignat}{5}\label{2_stage_5_order_ETSGLM_cond}
  -(1-v)           &+\widetilde{b}_1     & &+\widetilde{b}_2              & &+b_1    & &+b_2       & &=\alpha,            \notag\\
  \frac{1-v}{2}    &-\widetilde{b}_1     & &-(1-c_2)\widetilde{b}_2       & &        & &+c_2   b_2 & &=\frac{\alpha^2}{2},\notag\\
  -\frac{1-v}{3}   &+\widetilde{b}_1     & &+(1-c_2)^2\widetilde{b}_2     & &        & &+c_2^2 b_2 & &=\frac{\alpha^3}{3},\notag\\
  \frac{1-v}{4}    &-\widetilde{b}_1     & &-(1-c_2)^3\widetilde{b}_2     & &        & &+c_2^3 b_2 & &=\frac{\alpha^4}{4},\notag\\
  -\frac{1-v}{5}   &+\widetilde{b}_1     & &+(1-c_2)^4\widetilde{b}_2     & &        & &+c_2^4 b_2 & &=\frac{\alpha^5}{5},\\
  -(1-u_2)         &+\widetilde{a}_{21}  & &+\widetilde{a}_{22}           & &+a_{21} & &           & &=\alpha,            \notag\\
  \frac{1-u_2}{2}  &-\widetilde{a}_{21}  & &-(1-c_2)\widetilde{a}_{22}    & &        & &           & &=\frac{\alpha^2}{2},\notag\\
  -\frac{1-u_2}{3} &+\widetilde{a}_{21}  & &+(1-c_2)^2\widetilde{a}_{22}  & &        & &           & &=\frac{\alpha^3}{3},\notag\\
  \frac{1-u_2}{4}  &-\widetilde{a}_{21}  & &-(1-c_2)^3\widetilde{a}_{22}  & &        & &           & &=\frac{\alpha^4}{4},\notag
\end{alignat}
where $\alpha\in[0,1]$ in the first five equations~\eqref{2_stage_5_order_ETSGLM_cond} and $\alpha\in[0,c_2]$ in other ones.

\noindent
The coefficients are defined by
\mathindent=0em
\begin{alignat}{5}\label{2_stage_5_order_ETSGLM_coeff}
               u_2&=\left( \alpha+1 \right) ^{2} \left( 1-2\,\alpha+{\frac {3{\alpha}^{2}}{2\,c_{{2}}-1}} \right),  \notag\\
                 v&=-{\frac { \left( \alpha+1 \right) ^{2} \left(  \left( 10\,\alpha-5 \right) {c_{{2}}}^{2}-15\,c_{{2}}{\alpha}^{2}+ \left( \alpha+1 \right)  \left( 6\,{\alpha}^{2}-3\,\alpha+1 \right)  \right) }{5\,{c_{{2}}}^{2}-1}},  \notag\\
\widetilde{a}_{21}&={\alpha}^{2} \left( \alpha+1 \right) -{\frac {{\alpha}^{2}\left( \alpha+1 \right) ^{2} \left( 3\,c_{{2}}-1 \right) }{2\,c_{{2}}\left( 2\,c_{{2}}-1 \right) }},  \notag\\
\widetilde{a}_{22}&={\frac {{\alpha}^{2} \left( \alpha+1 \right) ^{2}}{2\,c_{{2}} \left( c_{{2}}-1 \right)  \left( 2\,c_{{2}}-1 \right) }},  \notag\\
   \widetilde{b}_1&={\frac {{\alpha}^{2} \left( \alpha+1 \right)  \left( 20{c_{{2}}}^{4}\!-\left(30 \alpha+10 \right) {c_{{2}}}^{3}+ \left( 12{\alpha}^{2}\!+\!3\alpha\!-\!13 \right) {c_{{2}}}^{2}+ \left( 4{\alpha}^{2}\!+\!11\alpha\!+\!3 \right) c_{{2}}-2\alpha\left( \alpha+1 \right)  \right)}{ 4c_2\left( 5{c_2}^{2}-1\right)  \left( c_2+1 \right) }},  \notag\\
   \widetilde{b}_2&={\frac { {\alpha}^{2}\left( \alpha+1 \right)^{2}\left( 5\,{c_{{2}}}^{2}- \left( 4\,\alpha-3 \right) c_{{2}}-2\,\alpha \right) }{4 c_2 \left( 5\,{c_{{2}}}^{2}-1 \right)\left( c_2-1 \right) }}\\
            a_{21}&=\alpha \left( \alpha+1 \right) ^{2} \left( 1-{\frac {\alpha\, \left( 3\,c_{{2}}-2 \right) }{2 \left( 2\,c_{{2}}-1 \right)  \left( c_{{2}}-1 \right) }} \right),  \notag\\
               b_1&={\frac {\alpha\left( \alpha+1 \right)^2  \left( 20{c_2}^{4}\!-\left(30 \alpha+20 \right) {c_{{2}}}^{3}+ \left( 12{\alpha}^{2}\!+\!21\alpha\!-\!4 \right) {c_{{2}}}^{2}+ \left( -4{\alpha}^{2}\!+\!3\alpha\!+\!4 \right) c_{{2}}-2\alpha\left( \alpha+1 \right)  \right)}{ 4c_2\left( 5{c_2}^{2}-1\right)  \left( c_2-1 \right) }},  \notag\\
               b_2&=-{\frac { {\alpha}^{2}\left( \alpha+1 \right)^{2}\left( 5\,{c_{{2}}}^{2}- \left( 4\,\alpha+7 \right) c_{{2}}+2\,\alpha+2 \right) }{4 c_2 \left( 5\,{c_{{2}}}^{2}-1 \right)\left( c_2+1 \right) }}.  \notag
\end{alignat}
To attain the discrete stage order five, we determine $c_2$ from $\Gamma_{2\,5}(1)=0$. We have\\
\begin{equation}
  \ c_2=\dfrac{11-\sqrt{41}}{10}.
\end{equation}
\mathindent=2.5em

\noindent
The relation $\Gamma_6(1)=-{\dfrac {16\ \ \bigl(17-2\sqrt {41}\,\bigr)}{75\,\bigl(71-11\sqrt {41}\,\bigr)}\ne 0}$
implies that it is impossible to attain discrete order six.\\

So we construct explicit TSGLM of uniform order five, uniform stage order four and discrete stage order five.

\begin{MyRemark}
  There is not a method of uniform stage order two in a class of explicit one-step RK methods for RFDEs.
  Indeed, for explicit one-step RK methods $c_1=0,\ c_2\ne0$ and  $a_{2,j}=0,\ j=2,\dots,s$,
  hence $\Gamma_{2\,k}=-\dfrac{\alpha^k}{k!}\ne0,\quad \alpha\in(0,c_2], \ k=2,3,\dots\ \ .$

  It is known~\cite{Dekker} that methods with low stage order suffer from the order reduction phenomenon when applied to stiff ODEs.
  Hence, explicit TSGLMs may be more appropriate for some mildly stiff RFDEs
  (of course, if the  smoothness conditions in Remark~\ref{Why_two_step_methods} and in the begining of section~\ref{sec_Order cond} hold).
\end{MyRemark}

\newpage
\renewcommand{\refname}{References}

\end{document}